\documentclass{article}[12pt]
\usepackage{graphicx} 
\usepackage[utf8]{inputenc}
\usepackage{authblk}
\usepackage{setspace}
\graphicspath{ {./figures/} }
\usepackage{subcaption}
\usepackage{amssymb,amsthm,amsmath}

\usepackage{color,xcolor}

\usepackage{hyperref}
\usepackage{lineno}

\usepackage{fullpage}

\newtheorem{theorem}{\bf Theorem}[section]

\newtheorem{lemma}[theorem]{\bf Lemma}
\newtheorem{corollary}[theorem]{\bf Corollary}

\newtheorem{conjecture}{\bf Conjecture}

\theoremstyle{plain}

\title{Rainbow matchings in non-uniform hypergraph systems \vspace{0.5 cm}}

\author{Ghaffar Raeisi\thanks{This work is based upon research funded by Iran National Science Foundation (INSF) under project No. 40400962 and also supported in part by a grant from IPM No. 1404050317.}\vspace{0.3 cm}\\
{\small
Department of Mathematical Sciences, Shahrekord University, Shahrekord, P.O.Box 115, Iran\\\vspace{0.2 cm}School of Mathematics, Institute for Research in Fundamental Sciences (IPM), P.O. Box: 19395--5746, Tehran, Iran\\ Email: \texttt{g.raeisi@sku.ac.ir}}
}

\date{}
\begin{document}
\maketitle
\begin{abstract}
Let $s\geq 1$, $r_1, r_2,\dots, r_{s+1}$ be given positive integers and $r=\max\{r_1, r_2,\ldots, r_{s+1}\}$. We prove that if $\mathcal{F} = \{\mathcal{F}_1, \ldots, \mathcal{F}_{s+1}\}$ is a family of hypergraphs on a ground set of $n$ vertices, $n\geq 4r^2(s+1)$, such that $\mathcal{F}_i\subseteq \binom{n}{r_i}$ and $|\mathcal{F}_i|>\binom{n}{r_i}-\binom{n-s}{r_i}-\binom{n-s-r}{r_i-1}+1$, for all $i \in [s+1]$, then $\mathcal{F}$ admits a rainbow matching of size $s+1$ or there exists $W \in \binom{[n]}{s}$ such that $W$ is a vertex cover of $\mathcal{F}_i$, for all $i \in [s+1]$. This improves upon the previously known range $n>2r^3s$  to $n>4r^2(s+1)$. Moreover, for a given family $\mathcal{F} = \{\mathcal{F}_1, \ldots, \mathcal{F}_{s+1}\}$ of hypergraphs $\mathcal{F}_i\subseteq \binom{n}{r_i}$, we give a best possible Ore-degree condition for $\mathcal{F}_{i}$'s that guarantees the existence of a rainbow matching of size $s+1$ for $\mathcal{F}$.
\end{abstract}

{\bf Key words:} Intersecting family, Rainbow matching, Ore-degree condition
\medskip
\section{ Introduction}

\medskip
An {\it $r$-uniform hypergraph} is a pair $(V,E)$, where $V$ is a finite set of vertices and $E\subseteq {V \choose r}$ is a family of $r$-subsets of
 $V$ called {\it edges} of ${\cal H}$. For a given hypergraph ${\cal H}$, we usually suppose that the vertex set of ${\cal H}$, $V({\cal H})$, is $[n]:=\{1,2,\dots, n\}$ and for simplicity, we identify $\mathcal{H}$ with its set of edges, and so by $e \in \mathcal{H}$ we mean that $e$ is an edge of $\mathcal{H}$. Thus $|\mathcal{H}|$ is the number of edges in $\mathcal{H}$. We often write $\mathcal{H} \subseteq \binom{[n]}{r}$ to define an $r$-uniform hypergraph $\mathcal{H}$ on vertex set $[n]$.
The {\it degree} of a vertex $v\in V({\cal H})$ is defined as the number of distinct edges containing $v$ and the {\it minimum degree} of ${\cal H}$ (resp.\ {\it maximum degree} of ${\cal H}$), $\delta({\cal H})$ (resp. $\Delta({\cal H})$), is the minimum (resp. the maximum) of $\deg(v)$ over all vertices $v\in V({\cal H})$. A {\it matching} in a hypergraph ${\cal H}$ is a set of pairwise disjoint edges and the {\it matching number} of ${\cal H}$, $\nu({\cal H})$, is the size of the largest matching in ${\cal H}$. A {\it vertex cover} of ${\cal H}$ is a set of vertices of ${\cal H}$ that meets all edges of ${\cal H}$, and the \emph{covering number} of ${\mathcal{H}}$, $\tau(\mathcal{H})$, is the minimum size of a vertex cover in $\mathcal{H}$. Clearly, for a given hypergraph ${\cal H}$, $\nu({\cal H})\leq \tau(\mathcal{H})$.

 \medskip
An $r$-uniform hypergraph ${\cal H}$ is called {\it intersecting}, if every pair of its edges has a non-empty intersection.
The celebrated theorem of Erd\H{o}s, Ko and Rado~\cite{EKR} states that when $n\geq 2r$, every intersecting $r$-uniform hypergraph on $[n]$ has at most ${n-1 \choose r-1}$ edges and moreover, for $n >2r$, the extremal hypergraph is unique up to isomorphism and it consists of all the $r$-subsets of $[n]$ that contain a fixed vertex $x\in [n]$.
We call such hypergraph a {\it 1-star}.

\medskip
In \cite{HM} Hilton and Milner proved a strong stability extension of the Erd\H{o}s--Ko--Rado theorem. An intersecting hypergraph ${\cal H}$ is  {\it trivial} if it is a subhypergraph of a $1$-star.  The {\it Hilton--Milner theorem} states that for $ n>2r$, every non-trivial intersecting $r$-uniform hypergraph on $n$ vertices has at most $\binom{n-1}{r-1} - \binom{n-r-1}{r-1} + 1$ edges. The corresponding extremal hypergraph ${\cal HM}_{n,r}$ consists of a fixed $r$-set $S$ and all $r$-subsets of $[n]$ that containing a fixed element $x \notin S$ and at least one element of $S$. In other words, for a given fixed $r$-set $S$,
\[
{\cal HM}_{n,r}=\Big\{H\in \binom{[n]}{r}:\quad x\in H,~ H\cap S\neq \emptyset\Big\}\cup S.
\]
In \cite{frk-kup}, the authors proved the following far-reaching generalization of the Hilton-Milner theorem.

\medskip
\begin{theorem}[Frankl and Kupavskii~\cite{frk-kup}]\label{franklkupavskii}
Suppose that $s\geq 1$, $r \ge 3$ and either $n \geq (s+\max\{25, 2s + 2\})r$ or $n \geq (2+o(1))sr$, where $o(1)$ is with respect to $s \to \infty$. Then, for any
$\mathcal{H} \subset \binom{[n]}{r}$ with $\nu(\mathcal{H}) = s < \tau(\mathcal{H})$ we have
\[
|\mathcal{H}| \le \binom{n}{r} - \binom{n - s}{r}-\binom{n - s - r}{r - 1}+1.
\]
\end{theorem}

\medskip
In general, the problem of finding a maximum matching in a hypergraph is an important problem in extremal combinatorics and has many applications in various different areas of mathematics and  computer science. Although the graph matching problem is fairly well-understood, and solvable in polynomial time, many problems related to hypergraph matching appear to be  difficult and remain unsolved. One of the most basic open questions in this area was raised in 1965 by Erd\H{o}s \cite{erdos}, who asked to determine the maximum possible number of edges that can appear in any $r$-uniform hypergraph ${\cal H}$ with matching number less than $s$. He conjectured that extremal constructions for this problem are one of the following two hypergraphs. The first is a clique consisting of all $r$-subsets of $rs-1$ vertices, which obviously has matching number $s-1$. The second one is $\mathcal{A}_0(n,r,s-1)$ which is an $r$-uniform hypergraph on $n$ vertices containing all edges that intersect a fixed set of $s-1$ vertices, which also forces the matching number to be $s-1$.

\begin{conjecture}[Erd\H os~\cite{erdos}]\label{conj}

If $n\geq rs-1$ and ${\cal H}$ is an $n$-vertex $r$-uniform hypergraph such that
\[
|{\cal H}|>\max\left\{{rs-1 \choose r}, {n \choose r}-{n-(s-1) \choose r}\right\},
\]
then $\nu({\cal H})\geq s$.
\end{conjecture}

For recent results about this conjecture, we refer the reader to references \cite{BDE, F2, F, FK2022, frku2017,FRR}. The rainbow version of the Erd\H os Matching Conjecture and the Hilton-Milner theorem were also widely studied by several authors. Let $\mathcal{F} = \{\mathcal{F}_1, \mathcal{F}_2, \ldots, \mathcal{F}_s\}$ be a family of given hypergraphs. A set of pairwise disjoint edges, one from each $\mathcal{F}_i$ is called a {\it rainbow matching}, and we say that $\mathcal{F}$ admits a rainbow matching of size $s$ if $\{\mathcal{F}_1, \ldots, \mathcal{F}_s\}$ has a rainbow matching.   The rainbow version of the Erd\H os Matching Conjecture was first studied by Huang, Loh and Sudakov, where they proved the following.

\medskip
\begin{theorem}[Huang, Loh and Sudakov \cite{HLS}]\label{HLS}
Let $n, r, s$ be positive integers such that $n > 3r^2 s$. If $\mathcal{F} = \{\mathcal{F}_1, \mathcal{F}_2, \ldots, \mathcal{F}_s\}$ is a family of $r$-uniform hypergraphs and
\[
|\mathcal{F}_i| > \binom{n}{r} - \binom{n-s+1}{r}
\]
for all $1 \le i \le s$, then $\mathcal{F}$ admits a rainbow matching.
\end{theorem}

\medskip
The bound $3r^2s$ on $n$ was improved to $12rs\log(e^2s)$ by Frankl and Kupavskii \cite{franku2020} and Lu,
Wang and Yu \cite{LWY} further improved this bound to $n > 2rs$, for $s$ sufficiently large. In \cite{LWY1}, the authors studied the rainbow version of Theorem \ref{franklkupavskii} and proved that if $s\geq 1$, $n\geq 2r^3s$ and $\mathcal{F} = \{\mathcal{F}_1, \mathcal{F}_2,\dots, \mathcal{F}_{s+1}\}$ is a family of $r$- uniform hypergraphs such that $|\mathcal{F}_i|>\binom{n}{r}-\binom{n-s}{r}-\binom{n-s-r}{r-1}+1$, for all $i \in [s+1]$, then $\mathcal{F}$ admits a rainbow matching of size $s+1$ or there exists $W \in \binom{[n]}{s}$ such that $W$ is a vertex cover of $\mathcal{F}_i$ for all $i \in [s+1]$. In fact, they proved the following more general result by considering non-uniform hypergraphs.

\medskip
\begin{theorem}[Lu, Wang and Yu~\cite{LWY1}]\label{rainbow}
Let $s\geq 1$ and $r_1, r_2,\ldots, r_{s+1}$  be given positive integers and $r=\max\{r_1, r_2,\ldots, r_{s+1}\}$. If $n\geq 2r^3s$ and $\mathcal{F} = \{\mathcal{F}_1, \dots, \mathcal{F}_{s+1}\}$ is a family of hypergraphs such that $\mathcal{F}_i\subseteq \binom{n}{r_i}$ and $|\mathcal{F}_i|>\binom{n}{r_i}-\binom{n-s}{r_i}-\binom{n-s-r}{r_i-1}+1$,
for all $i \in [s+1]$, then $\mathcal{F}$ admits a rainbow matching of size $s+1$ or there exists $W \in \binom{[n]}{s}$ such that $W$ is a vertex cover of $\mathcal{F}_i$, for all $i \in [s+1]$.
\end{theorem}

\medskip
In \cite{BCR}, the authors considered the rainbow version of Conjecture \ref{conj} via an Ore-type degree condition.  For an $r$-set $S \subseteq V({\cal H})$,
the {\it degree of} $S$ is defined as $\deg(S)=\sum_{v \in S}\deg(v)$ and the minimum of $\deg(S)$ over all non-edge $r$-subsets $S \not \in E(\mathcal{H})$ of $V({\cal H})$ is the \emph{Ore-degree} of $\mathcal{H}$ and is denoted by ${\sigma_r}({\cal H})$. For an ordinary graph $H$, we recover the standard definition of Ore degree:  $\sigma_2(H)=\min_{uv\notin E(H)}(\deg(u)+\deg(v))$. Clearly, for the family $\mathcal{A}_0(n,r,s)$ we have $\sigma_{r}(\mathcal{A}_0(n,r,s))=r(\binom{n-1}{r-1}-\binom{n-s-1}{r-1}).$ In \cite{BCR}, the authors considered an Ore-degree version of Theorem~\ref{HLS} and proved the following Ore-degree condition for the existence of a rainbow matching.

\medskip
 \begin{theorem}[Balogh, Palmer and Raeisi~\cite{BCR}]\label{ore}
Let $r\geq 3$, $s\geq 1$ and $n>3r^2(s+1)$. If $\mathcal{F} = \{\mathcal{F}_1, \dots, \mathcal{F}_{s+1}\}$ is a family of $r$-uniform hypergraphs on $n$ vertices such that $\sigma_r({\cal F}_i)> r\left({n-1 \choose r-1}-{n-s-1 \choose r-1}\right)$, for every $i\in [s+1]$, then ${\cal F}$ admits a rainbow matching.
\end{theorem}

\medskip
Note that the Ore-degree condition in Theorem \ref{ore} is best possible,  because of $\mathcal{A}_0(n,r,s)$. In this paper, we first improve the range $n\geq 2r^3s$ in Theorem \ref{rainbow} to $n\geq 4r^2(s+1)$ and also we extend Theorem \ref{ore} to families of non-uniform hypergraphs. In other words, the following are the main results of the paper.

\medskip
\begin{theorem}\label{rainbowhilton}
Let $s\geq 1$ and $r_1, r_2,\ldots, r_{s+1},n$ be given positive integers and $r=\max\{r_1, r_2,\ldots, r_{s+1}\}$. If $\mathcal{F} = \{\mathcal{F}_1, \dots, \mathcal{F}_{s+1}\}$ is a family of hypergraphs on a ground set of $n$ vertices, $n\geq 4r^2(s+1)$, such that $\mathcal{F}_i\subseteq \binom{n}{r_i}$ and $|\mathcal{F}_i|>\binom{n}{r_i}-\binom{n-s}{r_i}-\binom{n-s-r}{r_i-1}+1$, for all $i \in [s+1]$, then either $\mathcal{F}$ admits a rainbow matching of size $s+1$ or there exists $W \in \binom{[n]}{s}$ such that $W$ is a vertex cover of $\mathcal{F}_i$, for all $i \in [s+1]$.
\end{theorem}

\medskip
It is worth noticing that the following family shows that the bound on $|\mathcal{F}_i|$ in Theorem \ref{rainbowhilton} is best possible. For given positive integers  $r, n , s$, let $\mathcal{H}^{(r)}(n, s)$ be the following family:
\[
\mathcal{H}^{(r)}(n, s) =
\left\{
H \in \binom{[n]}{r} : H \cap [s-1] \neq \varnothing
\right\}
\cup
\left\{
H \in \binom{[s, n]}{r} : s \in H,\; H \cap [s+1, s+r] \neq \varnothing
\right\}
\cup  [s+1, s+r].
\]
Here, for $s<n$, by $[s,n]$ we mean the set $\{s,s+1,\ldots, n\}$. The family $\mathcal{H}^{(r)}(n, s)$ is an $r$-uniform hypergraph which is the union of $s-1$ stars with centers in $\{1,2,\ldots, s-1\}$ and the Hilton--Milner family on $\{s,s+1,\ldots, n\}$. Clearly, $\nu(\mathcal{H}^{(r)}(n, s))\leq s$ and $|\mathcal{H}^{(r)}(n, s)|=\binom{n}{r}-\binom{n-s}{r}-\binom{n-s-r}{r-1}+1$. Choosing $\mathcal{F}_1= \cdots= \mathcal{F}_{s+1}=\mathcal{H}^{(r)}(n, s)$, shows that the lower bound on $|\mathcal{F}_i|$ in Theorem \ref{rainbowhilton} is best possible.

\medskip
\begin{theorem}\label{main 2}
Let $s\geq 1$ and  $r_1, r_2,\ldots, r_{s+1}$ be given positive integers, $r=\max\{r_1, r_2,\dots, r_{s+1}\}$, $\min_i r_i>s$ and $n\geq 4r^2(s+1)$. If $\mathcal{F} = \{\mathcal{F}_1, \dots, \mathcal{F}_{s+1}\}$ be is a family of hypergraphs such that $\mathcal{F}_i\subseteq \binom{n}{r_i}$ and $\mathcal{F}$ does not admit a rainbow matching of size $s+1$, then $$\displaystyle\min_{i}\sigma_{r_i}(\mathcal{F}_i)\leq r_i\left(\binom{n-1}{r_i-1}-\binom{n-s-1}{r_i-1}\right).$$
\end{theorem}

As an immediate and easy corollary of Theorem \ref{main 2}, we have the following corollary on cross-intersecting families. Two families ${\cal A}$ and ${\cal B}$  are called {\it cross-intersecting} if for every $A \in {\cal A}$  and $B \in {\cal B}$, their intersection $A\cap B$ is non-empty. Noting that for $\mathcal{A}=\mathcal{B}$ the cross-intersecting property reduces to $\mathcal{A}$ being intersecting.

\medskip
\begin{corollary}
For given $r,k\geq 2$, $l=\max\{r,k\}$ and $n\geq 8l^2$, if ${\cal A}\subseteq \binom{[n]}{r}$ and ${\cal B} \subseteq \binom{[n]}{k}$ are cross-intersecting, then $\sigma_r({\cal A})\leq r{n-2 \choose r-2}$ or $\sigma_k({\cal B})\leq k{n-2 \choose k-2}$. In particular, $\delta({\cal A})\leq {n-2 \choose r-2}$ or $\delta({\cal B})\leq {n-2 \choose k-2}$.
\end{corollary}
%
\bigskip
\section{Proofs of the main results}

In this section we prove the main results of the paper. First, we list some lemmas and theorems that we use in our proofs. The following result of M\"{o}rs extending the Hilton-Milner Theorem to two cross-intersecting families.

\begin{theorem}[M\"{o}rs~\cite{Mors}]\label{Mors}
Let $n, r, k$ be positive integers such that $n \ge r + k$.
Let $\mathcal{A} \subseteq \binom{[n]}{r}$ and $\mathcal{B} \subseteq \binom{[n]}{k}$.
Suppose $\bigcap_{C \in \mathcal{A} \cup \mathcal{B}} C = \emptyset$ and for all $A \in \mathcal{A}, \ B \in \mathcal{B}$, $A \cap B \neq \emptyset$.
Then either
\[
|\mathcal{A}| \leq \binom{n-1}{r-1}-\binom{n-1-k}{r-1}+1,
\]
or
\[
|\mathcal{B}| \leq \binom{n-1}{k-1}-\binom{n-1-r}{k-1}+1.
\]
\end{theorem}

In \cite{BCR}, the authors proved the following, giving a tight lower bound on the size of an $n$-vertex $r$-uniform hypergraph ${\cal  H}$ in terms of $\sigma_r({\cal H})$. As this lemma is used frequently in this paper, for the sake of completeness, a proof is given here.

\medskip
\begin{lemma}\label{lem}
Let ${\cal  H} $ be an $r$-uniform hypergraph on $n$ vertices. Then  $|{\cal H}|\geq \frac{n}{r^2}\cdot \sigma_r({\cal H})$, with equality only if $\mathcal{H}$ is a  regular hypergraph.
\end{lemma}

\begin{proof}
For each of the $\binom{n}{r}-|{\cal H}|$ non-edge $r$-sets $S \notin {\cal H}$, we have $\sum_{v \in S} \deg(v) \geq \sigma_r({\cal H}).$ Summing over all non-edge $r$-sets, we obtain
   \[
   \sum_{S \in  \binom{[n]}{r} \setminus {\cal H}} \sum_{v \in S} \deg(v) \geq \sigma_r({\cal H})  \left(\binom{n}{r}-|{\cal H}|\right).
   \]
Since the number of non-edge $r$-sets containing a vertex $v$ is $\binom{n-1}{r-1}-\deg(v)$, the left-hand side counts each vertex degree $\deg(v)$ exactly $\binom{n-1}{r-1}-\deg(v)$ times. Thus,
   \[
   \sum_{v \in V({\cal H})} \deg(v) \left(\binom{n-1}{r-1}-\deg(v)\right) \geq \sigma_r({\cal H}) \left(\binom{n}{r}-|{\cal H}|\right).
   \]
Now, $\sum_{v \in V(\mathcal{H})} \deg(v)= r|{\cal H}|$, hence
   \[
   r |{\cal H}| \binom{n-1}{r-1}-\sum_{v\in V(\mathcal{H})} \deg(v)^2 \geq \sigma_r({\cal H}) \left(\binom{n}{r}-|{\cal H}|\right).
   \]
Applying the Cauchy-Schwartz inequality gives $$\sum \deg(v)^2 \geq \frac{1}{n}\left(\sum \deg(v)\right)^2=\frac{r^2 |{\cal H}|^2}{n} .$$ Therefore,

  \[
  r|{\cal H}| \binom{n-1}{r-1}-\frac{r^2 |{\cal H}|^2}{n} \geq \sigma_r({\cal H})  \left(\binom{n}{r}-|{\cal H}|\right),
  \]
and equivalently,
\[
0\geq \frac{r^2}{n}|{\cal H}|^2-\left(\sigma_r({\cal H})+r\binom{n-1}{r-1}\right)|{\cal H}|+\sigma_r({\cal H})\binom{n}{r}
 =\left(\frac{r^2}{n}|{\cal H}|-\sigma_r(\mathcal{H})\right)\left(|{\cal H}|-\binom{n}{r}\right).
   \]
Thus, $|{\cal H}|\geq \frac{n}{r^2}\cdot \sigma_r({\cal H})$. Clearly, equality holds only if all vertex degrees of ${\cal H}$ are equal.
\end{proof}

\medskip
In \cite{LWY1}, the authors used the Kruskal-Katona theorem to prove the following, which give a connection between the number of edges in an $r$-uniform hypergraph and its $k$-uniform 'expansion'.

\medskip
\begin{lemma}[Lu, Wang and Yu~\cite{LWY1}]\label{lema0}
Let $n\geq k\geq r\geq 2$ be given. Also, let $\mathcal{F}\subseteq \binom{n}{r}$ and $\mathcal{G}= \{f\cup e: f\in \mathcal{F}, e\in \binom{[n]\setminus f}{k-r}\}$. If
$|\mathcal{F}|>\binom{n}{r}-\binom{n-s}{r}-\binom{n-s-k}{r-1}+1$, then $|\mathcal{G}|>\binom{n}{k}-\binom{n-s}{k}-\binom{n-s-k}{k-1}+1$.
\end{lemma}

\medskip
Before giving a proof for Theorem \ref{rainbowhilton}, we prove the following lemma.

\medskip
\begin{lemma}\label{lema}
For given $r,s\geq 2$ and $n\geq 4r(r-1)(s+1)+s$,
\[
\binom{n}{r}-\binom{n-s}{r}-\binom{n-s-r}{r-1}+1>s^{2}(3r-2)\binom{n-2}{r-2}.
\]
\end{lemma}
\begin{proof}
First, note that
\[
\binom{n}{r}-\binom{n-s}{r}=\sum_{i=0}^{s-1} \binom{n-1-i}{r-1}\geq s\binom{n-s}{r-1}.
\]
Also, $\binom{n-s-r}{r-1} \leq \binom{n-s}{r-1}$ and so
\begin{equation*}\label{1}
\binom{n}{r}-\binom{n-s}{r}-\binom{n-s-r}{r-1}+1>(s-1)\binom{n-s}{r-1}.
\end{equation*}
Therefore, it is sufficient to prove that $(s-1)\binom{n-s}{r-1}\geq s^{2}(3r-2)\binom{n-2}{r-2}$, or equivalently, $\frac{(s-1)\binom{n-s}{r-1}}{\binom{n-2}{r-2}}>s^{2}(3r-2)$. Clearly,
\begin{equation}\label{2}
\frac{\binom{n-s}{r-1}}{\binom{n-2}{r-2}}=\frac{n-s}{r-1}.\frac{\binom{n-s-1}{r-2}}{\binom{n-2}{r-2}}=\frac{n-s}{r-1}\prod_{j=1}^{r-2} \frac{n-s-j}{n-1-j}=\frac{n-s}{r-1}\prod_{j=1}^{r-2}(1-\frac{s-1}{n-1-j})\geq \frac{n-s}{r-1}(1-\frac{s-1}{n-r+1})^{r-2}.
\end{equation}
Since, $(1-x)^m>1-mx$, for $0<x<1$, and $n-s>4r(r-1)(s+1)$, from \eqref{2} we obtain that
\begin{equation*}
\frac{\binom{n-s}{r-1}}{\binom{n-2}{r-2}}> \frac{n-s}{r-1}\left(1-\frac{(r-2)(s-1)}{n-r+1}\right)>4r(s+1)\left(1-\frac{1}{4r}\right),
\end{equation*}
where the last inequality follows from $n-(r-1)> 4r(r-2)(s-1).$ Therefore,
\begin{equation*}
\frac{(s-1)\binom{n-s}{r-1}}{\binom{n-2}{r-2}}> (s^2-1)(4r-1)>s^{2}(3r-2),
\end{equation*}
provided that $r,s\geq 2$. This observation completes the proof of theorem.
\end{proof}

\medskip
To prove Theorem \ref{rainbowhilton}, we first prove an analogous  theorem for uniform hypergraphs.

\medskip
\begin{theorem}\label{rainbowhiltonuniform}
Let $s, r$ be given positive integers and $n\geq 4r^2(s+1)$. If $\mathcal{F} = \{\mathcal{F}_1, \ldots, \mathcal{F}_{s+1}\}$ is a family of $r$-uniform hypergraphs on a ground set of $n$ vertices and $|\mathcal{F}_i|>\binom{n}{r}-\binom{n-s}{r}-\binom{n-s-r}{r-1}+1$, for all $i \in [s+1]$, then either $\mathcal{F}$ admits a rainbow matching of size $s+1$ or there exists $W \in \binom{[n]}{s}$ such that $W$ is a vertex cover of $\mathcal{F}_i$, for all $i \in [s+1]$.
\end{theorem}
\begin{proof}
We use induction on $s$ to prove the theorem. The base case $s=1$ follows from Theorem \ref{Mors} and thus, we may assume that $s\geq 2$. We suppose that
$\mathcal{F}$ does not admit a rainbow matching and then we prove that there exists $W \in \binom{[n]}{s}$ such that $W$ is a vertex cover of $\mathcal{F}_i$, for all $i \in [s+1]$. First, we prove that for each $i$, $1\leq i\leq s+1$, we have  $\Delta({\cal F}_i)\leq rs{n-2 \choose r-2}$. Indeed, suppose that some ${\cal F}_{i}$, say ${\cal F}_{s+1}$, has a vertex $v$ of degree greater than $rs{n-2 \choose r-2}$. For $i\in [s]$, define ${\cal F}_i'$ as the subhypergraph of ${\cal F}_{i}$ induced by the vertex set $V({\cal F}_{i}) \setminus \{v\}$.
Clearly, for $i\in [s]$,
\[
|\mathcal{F}_i'|\geq |\mathcal{F}_i|-\binom{n-1}{r-1} > \binom{n-1}{r} - \binom{(n-1)-(s-1)}{r}-\binom{(n-1)-(s-1)-r}{r-1}+1,
\]
and so by the induction hypothesis applied to the family $\mathcal{F}'=\{\mathcal{F}_1', \dots, \mathcal{F}_s'\}$, we obtain that $\mathcal{F}'$ admits a rainbow matching or there exists
$W' \in \binom{[n]\setminus \{v\}}{s-1}$ such that $W'$ is a vertex cover of $\mathcal{F}_i'$ for all $i \in [s]$. If $\mathcal{F}'$ admits a rainbow matching $\mathcal {M}$, then $\deg_{\mathcal{F}_{s+1}}(v)>rs{n-2 \choose r-2}$ implies that there exists an edge $e\in \mathcal{F}_{s+1}$ such that $e\cap V(\mathcal {M})=\emptyset$. Now $\mathcal {M}\cup e$ is a rainbow matching for $\mathcal {F}$ and we are done. Therefore, there exists $W' \in \binom{[n]\setminus \{v\}}{s-1}$ such that $W'$ is a vertex cover of $\mathcal{F}_i'$, for all $i \in [s]$.

\medskip
Set $W = W' \cup \{v\}$. If every edge of $\mathcal{F}_{s+1}$ meets $W$, then $W$ is a vertex cover for each $\mathcal{F}_i$, $i\in [s+1]$, and we are done. Otherwise, choose $f \in \mathcal{F}_{s+1}$ with $f \cap W = \emptyset$. For each $i \in [s]$, using the inclusion-exclusion principle, the number of edges of $\mathcal{F}_i$ intersecting both $f$ and $W$, say $e_i(W,f)$, is at most
\[
\binom{n}{r}-\binom{n-r}{r}-\binom{n-s}{r}+\binom{n-s-r}{r}.
\]
Thus, the number of edges of $\mathcal{F}_i$, $i \in [s]$, disjoint from $f$ is
\begin{align}
|\mathcal{F}_i\setminus f| &> \binom{n}{r}-\binom{n-s}{r}-\binom{n-s-r}{r-1}+1-e_i(W,f)\nonumber\\
&= \binom{n-r}{r}-\left(\binom{n-r-s}{r-1}+\binom{n-r-s}{r}\right)+1\\
&> \binom{n-r}{r}-\binom{(n-r)-(s-1)}{r}.
\end{align}

Now, by Theorem \ref{HLS}, $\{\mathcal{F}_1\setminus f, \dots, \mathcal{F}_s\setminus f\}$ admits a rainbow matching $\mathcal{M}'$ of size $s$. Then $\mathcal{M}' \cup \{f\}$ is a rainbow matching of size $s+1$ for the family $\mathcal{F}$.

\medskip
Therefore, for each $i\in [s+1]$, $\Delta({\cal F}_i)\leq rs{n-2 \choose r-2}$. After re-indexing the vertices of ${\cal F}_i$, we may assume that $rs\binom{n-2}{r-2}\geq \deg(v_1)\geq \deg(v_2) \geq \cdots \geq \deg(v_n)$. Now, we distinguish two cases based on the $(s+1)$-th largest degree in the ${\cal F}_{i}$'s.

\medskip
\noindent
{\it{ Case 1:} For some $i$, $1\leq i\leq s+1,$ the $(s+1)$-th largest degree in ${\cal F}_{i}$ is at most $2s\binom{n-2}{r-2}$.}

\medskip
Let the $(s+1)$-th largest degree in ${\cal F}_{s+1}$ be at most $2s\binom{n-2}{r-2}$. Clearly, for each $i\in [s]$,
\begin{align}
|\mathcal{F}_i| &>\binom{n}{r}-\binom{n-s}{r}-\binom{n-s-r}{r-1}+1 \nonumber\\
&> \binom{n}{r}-\binom{n-s}{r}-\binom{n-s}{r-1}~~~~~~~~(\text{because}~ \binom{n-s}{r-1}\geq \binom{n-s-r}{r-1})\nonumber\\
&> \binom{n}{r}-\binom{n-(s-1)}{r}~~~~~~~~(\text{because}~\binom{n-(s-1)}{r}= \binom{n-s}{r}+\binom{n-s}{r-1}).\nonumber
\end{align}
and so by Theorem \ref{HLS}, $\{\mathcal{F}_1, \mathcal{F}_2,\ldots \mathcal{F}_{s}\}$ admits a rainbow matching, say $\mathcal{M}$ which spans $rs$ vertices. Since the $(s+1)$-th largest degree of ${\cal F}_{s+1}$ is at most $2s\binom{n-2}{r-2}$, the sum of degrees of the $rs$ vertices of $V(\mathcal{M})$ in ${\cal F}_{s+1}$ is at most
\[
s\left[rs\binom{n-2}{r-2}\right]+s(r-1)\left[2s\binom{n-2}{r-2}\right]=s^{2}(3r-2)\binom{n-2}{r-2}.
\]
Now, by Lemma \ref{lema},  $|{\cal  F}_{s+1}|> s^{2}(3r-2)\binom{n-2}{r-2}$, for $n>4r^2(s+1)$, means that in $\mathcal{F}_{s+1}$ there is an edge $e$ disjoint from $\mathcal{M}$. Therefore, $\mathcal{M}\cup e$ produces an $(s+1)$-rainbow matching for $\mathcal{F}=\{{\cal F}_{1}, \dots, {\cal F}_{s+1}\}$.

\medskip
\noindent
{\it { Case 2:} For all $i$, $1\leq i\leq s+1$, ${\cal  F}_{i}$ has at least $s+1$ vertices of degree greater than $2s\binom{n-2}{r-2}$.}

\medskip
For each $i\in [s+1]$, select distinct vertices $v_{i}\in V({\cal F}_{i})$ such that the degree of $v_{i}$ in ${\cal  F}_{i}$ is greater than $2s\binom{n-2}{r-2}$. For each $i$, consider all subsets of $V({\cal  F}_{i})\setminus\{v_{1},\ldots,v_{s+1}\}$ which together with $v_{i}$ form an edge of ${\cal  F}_{i}$. Denote the corresponding $(r-1)$-uniform hypergraph by ${\cal  F}'_{i}$.
For any fixed $1\leq i\leq s+1$ and $j\neq i$, there are at most $\binom{n-2}{r-2}$ edges of ${\cal  F}_i$ containing both $v_i$ and $v_j$. Therefore, for every hypergraph ${\cal  F}_i$, $i\in[s+1]$, we have
\[
|{\cal  F}^{\prime}_{i}|>2s\binom{n-2}{r-2}-s\binom{n-2}{r-2}=s\binom{n-2}{r-2}>\binom{n-s-1}{r-1}-\binom{n-2s-1}{r-1}.
\]
Now, applying Theorem~\ref{HLS} to the $(r-1)$-uniform hypergraphs $\mathcal{F}_1',\dots, \mathcal{F}_{s+1}'$ on $n-s-1$ vertices, we obtain an $(s+1)$-rainbow matching $e_{1} \in \mathcal{F}'_1, e_2 \in \mathcal{F}'_2, \dots, e_{s+1} \in \mathcal{F}'_{s+1}$. Taking the edges $e_{i}\cup\{v_{i}\}$ for all $i\in[s+1]$, we obtain an $(s+1)$-rainbow matching in ${\cal F}_{1},\dots, \mathcal{F}_{s+1}$, completing the proof of theorem.
\end{proof}

\medskip
Now, we are ready to give a proof for Theorem \ref{rainbowhilton}.

\medskip
{\it \bf Proof of Theorem \ref{rainbowhilton}:}
To prove the theorem, first we reduce the problem of finding rainbow matchings for non-uniform hypergraphs to the one for uniform hypergraphs. For each $i\in [s+1],$ set $\mathcal{G}_i= \{e\cup f: e\in \mathcal{F}_i, f\in \binom{[n]\setminus e}{r-r_i}\}$. Clearly, for each $i\in [s+1],$ $\mathcal{G}_i$ is an expansion of $\mathcal{F}_i$ and is an $r$-uniform hypergraph on $n$ vertices. By Lemma \ref{lema0}, since $|\mathcal{F}_i|>\binom{n}{r_i}-\binom{n-s}{r_i}-\binom{n-s-r}{r_i-1}+1$, we obtain that $|\mathcal{G}_i|>\binom{n}{r}-\binom{n-s}{r}-\binom{n-s-r}{r-1}+1$. Now, having applied Theorem \ref{rainbowhiltonuniform} to the family $\mathcal{G}=\{\mathcal{G}_1, \mathcal{G}_2,\ldots \mathcal{G}_{s+1}\}$, we obtain that $\mathcal{G}$ admits an $(s+1)$-rainbow matching or there exists $W \in \binom{[n]}{s}$ such that $W$ is a vertex cover of $\mathcal{G}_i$ for all $i \in [s+1]$. Note that the existence of a $(s+1)$-rainbow matching for the family $\mathcal{G}$ implies the existence of an $(s+1)$-rainbow matching for $\mathcal{F}$. Also if $W$ is a vertex cover of every $\mathcal{G}_i$, then certainly this would be vertex cover for each $\mathcal{F}_i$. To see this, on the contrary, let $e \in \mathcal{F}_i$ be an edge disjoint from $W$. Choose some $f \subseteq [n] \setminus (e \cup W)$ of size $r- r_i$ (such an $f$ exists because $n$ is large enough). Then $e \cup f$ is an edge of $\mathcal{G}_i$ disjoint from $W$, contradicting that $W$ is a vertex cover of $\mathcal{G}_i$. This observation completes the proof of theorem.
$\hfill \blacksquare$

\medskip
Now, we give a proof for Theorem \ref{main 2}.

\medskip
{\it \bf Proof of Theorem \ref{main 2}:}
On the contrary, let for all $i \in [s+1]$, $\sigma_{r_i}(\mathcal{F}_i)>r_i\big(\binom{n-1}{r_i-1}-\binom{n-s-1}{r_i-1}\big)$. To get a contradiction, first, we prove that for $n>4r^2(s+1)$, $r\geq r_i>s$, $|{\cal  F}_i|>\binom{n}{r_i} - \binom{n-s}{r_i} - \binom{n-s-r}{r_i-1}+1$, for each $i$, $1\leq i\leq s+1$. Since $\sigma_{r_i}(\mathcal{F}_i)>r_i\big(\binom{n-1}{r_i-1}-\binom{n- s-1}{r_i-1}\big)$, by Lemma \ref{lem}, we obtain that
\[
|\mathcal{F}_i|>\frac{n}{r_i}\left({n-1 \choose r_i-1}-{n-s-1 \choose r_i-1}\right).
\]

Therefore, it is sufficient to prove that for $n>4r^2(s+1)$, $r\geq r_i>s$,
\begin{align}
\frac{n}{r_i}\left({n-1 \choose r_i-1}-{n-s-1 \choose r_i-1}\right)>\binom{n}{r_i}-\binom{n-s}{r_i}-\binom{n-s-r}{r_i-1}+1.\label{0}
\end{align}

Set $A:={n-s-1 \choose r_i-1},$ $B={n-s-r \choose r_i-1}$. Now, using $\frac{n}{r_i}{n-1 \choose r_i-1}={n \choose r_i}$ and ${n-s \choose r_i}=\frac{n-s}{r_i}{n-s-1 \choose r_i-1}$, equation \eqref{0} is equivalent to

\begin{align}
B-\frac{s}{r_i}A-1>0.\label{00}
\end{align}

Clearly,

\begin{align}
\frac{B}{A} =\frac{{n-s-r \choose r_i-1}}{{n-s-1 \choose r_i-1}}= \prod_{j=0}^{r-2} \frac{n-s-r_i-j}{n-s-j-1}= \prod_{j=0}^{r-2} \left(1-\frac{r_i-1}{n-s-j-1}\right).\label{1}
\end{align}

Using $\prod_{j=1}^{m} (1 - x_j) \ge 1 - \sum_{j=1}^{m} x_j$, $0\leq x_j\leq 1$, from equation \eqref{1} we obtain that

\begin{align}
\frac{B}{A}\geq 1-\frac{(r_i-1)(r-1)}{n-s-r_i+1}>1-\frac{(r_i-1)(r-1)}{2r^2s}.\label{20}
\end{align}

(The last inequality follows from the fact that for $n\geq 4r^2s$ and $r>r_i$, $n-s-r_i+1>4r^2(s+1)-s-r+1> 2{r}^2s$). Now, $s(r_i-s)-(r_i-1)=(s-1)(r_i-s-1)\geq 0,$
or equivalently, $s(r_i-s)\geq (r_i-1)$. Since $r\geq r_i$, we obtain that $rs(r_i-s)\geq r_i(r_i-1)$ and thus, $r^2s(r_i-s)\geq r_i(r_i-1)(r-1)$. This means that, $\frac{(r_i-1)(r-1)}{2r^2s}<\frac{{r}_i-s}{2r_i}$ and so by \eqref{20}, we obtain
\[
\frac{B}{A}>1-\frac{{r}_i-s}{2r_i}=\frac{s}{r_i}+\frac{r_i-s}{2r_i}.
\]
Therefore,
\[
B-\frac{s}{r_i}A>\frac{(r_i-s)A}{2r_i}\geq \frac{A}{2r_i}.
\]
But, $A={n-s-1 \choose r_i-1}\geq n-s-1>4r^2s-s-1\geq4r_i^2-r_i-1>2r_i$, and so $B-\frac{s}{r_i}A>1$, proving \eqref{00}. Therefore, for each $i\in[s+1]$,
\[
|{\cal  F}_i|>\binom{n}{r_i}-\binom{n-s}{r_i}-\binom{n-s-r}{r_i-1}+1.
\]
Note that $\tau(\mathcal{F}_i) \ge s+1$, for each $i\in [s+1]$. Indeed, if $\tau(\mathcal{F}_i)\leq s$ for some $i$, $i\in [s+1]$, then $\mathcal{F}_i$ is isomorphic to a subfamily of $\mathcal{A}_0(n,r_i,s)$ and so $\sigma_{r_i}({\cal  F}_i)\leq \sigma_{r_i}(\mathcal{A}_0(n,r_i,s)) = r_i\left({n-1 \choose r_i-1}-{n-s-1 \choose r_i-1}\right)$, which contradicts the assumption. Therefore,  $\tau(\mathcal{F}_i) \ge s+1$, for each $i$, $i\in [s+1]$. Now, using Theorem \ref{rainbowhilton}, we conclude that the family $\mathcal{F}$ admits a rainbow matching of size $s+1$, a contradiction. This contradiction completes the proof of theorem.

$\hfill \blacksquare$

\bigskip
\noindent{\bf  Conflict of interest:} The author states that there is no conflict of interest.

\end{document}